\documentstyle[12pt,amssymb]{article}
\newcommand{\R}{{\Bbb R}}
\newcommand{\C}{{\Bbb C}}
\newcommand{\N}{{\Bbb N}}
\newcommand{\bB}{{\bf B}}

\newcommand{\bH}{{\bf H}}

\newcommand{\bR}{{\bf R}}

\newcommand{\cD}{{\cal D}}

\newcommand{\cH}{{\cal H}}
\newcommand{\cI}{{\cal I}}

\newcommand{\cN}{{\cal N}}

\newcommand{\cV}{{\cal V}}
\newcommand{\ri}{{\rm i}}

\newcommand{\sP}{{\sf P}}
\newcommand{\sQ}{{\sf Q}}
\newcommand{\sm}{{m}}

\newcommand{\Min}{\mathop{{\rm min}}\limits}
\newcommand{\Max}{\mathop{{\rm max}}\limits}
\newcommand{\Inf}{\mathop{\rm inf}}
\newcommand{\Sup}{\mathop{\rm sup}}

\newcommand{\reduction}[2]{#1 \biggr|_{#2}}

\newcommand{\Lim}{\mathop{{\rm lim}}\limits}
\newcommand{\diag}{\mathop{\rm diag}}
\newcommand{\dist}{\mathop{\rm dist}}

\newcommand{\Real}{\mathop{\rm Re}}
\newcommand{\Inter}{\mathop{\rm Int}}
\newcommand{\lal}{{\langle}}
\newcommand{\ral}{{\rangle}}

\newtheorem{theorem}{\sc Theorem}
\newtheorem{lemma}{\sc Lemma}

\newtheorem{remark}{\sc Remark}

\title
{\normalsize\bf
\vskip 2truecm
OPERATOR INTERPRETATION OF RESONANCES
GENERATED BY SOME OPERATOR MATRICES%
\footnote{LANL E-print {\tt math.SP/9809093}.
Contribution to Proceedings of the
Mark Krein International Conference on Operator Theory
and Applications, Odessa, August 18-22, 1997.}\,\,%
\footnote{Financial support of this work by the DFG, INTAS
and RFBR is kindly acknowledged.}
}
\author
{\normalsize
R. MENNICKEN\thanks{{\it Address:} Naturwiss.
      Fakult\"at I -- Mathematik, Universit\"at Regensburg,
      D-93040 Regensburg, Germany. {\it E-mail}:
      {\tt reinhard.mennicken@mathematik.uni-regensburg.de}.}\,,
A. K. MOTOVILOV\thanks{{\it Address}:
      Physikalishes Institut, Universit\"at Bonn,
      Endenicher Allee 11-13, D-53115 Bonn, Germany. \,\,
      {\it On leave of absence from}
      the Bogoliubov Laboratory of Theoretical Physics,
      Joint Institute for Nuclear Research,
      141980 Dubna, Russia. {\it E-mail}:
      {\tt motovilov@physik.uni-bonn.de}
      and {\tt motovilv@thsun1.jinr.ru}.}
}

\date{}
\begin{document}
\maketitle
\thispagestyle{empty}

\baselineskip=12pt
\begin{quote}
We consider the analytic continuation of the transfer function
for a $2\times2$ matrix Hamiltonian into the unphysical sheets
of the energy Riemann surface. We construct a family of
non-selfadjoint operators which reproduce certain parts of the
transfer-function spectrum including resonances situated on the
unphysical sheets neighboring the physical sheet. On this basis,
completeness and basis properties for the root vectors of the
transfer function (including those for the resonances) are
proved.
\end{quote}

\vskip 1truecm
\baselineskip=15pt
\section{Introduction}
\label{Intro}
In this paper we deal with $2\times2$ operator matrices
\begin{equation}
\label{twochannel}
\bH=\left(
\begin{array}{cc}
                     A_0              &          B_{01}          \\
                     B_{10}           &          A_1
\end{array}
\right)
\end{equation}
acting in an orthogonal sum  ${\cal H}={\cal H}_0\oplus{\cal H}_1$
of separable Hilbert spaces ${\cal H}_0$ and ${\cal H}_1$.  The
entries $A_0:{\cal H}_0\rightarrow{\cal H}_0$, and
$A_1:{\cal H}_1\rightarrow{\cal H}_1$, are assumed to be
self-adjoint operators with
domains ${\cal D}(A_0)$ and ${\cal D}(A_1)$, respectively.  It
is assumed that the couplings $B_{ij}:
\cH_j\rightarrow\cH_i$, $i,j=0,1$, $i{\neq}j$, are bounded
operators  (i.\,e., $B_{ij}\in\bB(\cH_j,\cH_i)$) and $B_{01}=B^*_{10}$.
Under these assumptions the matrix $\bH$ is a self-adjoint operator
in $\cH$ with domain $\cD(\bH)=\cD(A_0)\oplus\cD(A_1)$.
Note that operators of the form~(\ref{twochannel}) arise in many
of quantum-physical problems (see
e.\,g.,~\cite{JaffeLow}\,--\,\cite{PavlovShushkov}).

In the spectral theory of operator matrices~(\ref{twochannel}) an
important role is played by the {\em transfer functions}
\begin{equation}
\label{trf}
M_i(z)=A_i-z+V_i(z) \quad\mbox{where}\quad
V_i(z)=-B_{ij}R_j(z)B_{ji}, \quad i=0,1, \quad j\neq i.
\end{equation}
A particular role of the functions $M_i(z)$ can be
understood already from the fact that the resolvent of the
operator~(\ref{twochannel}) can be expressed explicitly in terms of
the inverse transfer functions $M_0^{-1}(z)$ or $M_1^{-1}(z)$.
Therefore, in studying the spectral properties of the transfer
functions one studies at the same time the spectral
properties of the operator matrix ${\bf H}$.

In the papers~\cite{McKellarMcCay}, \cite{SchmidEW} the following
question was raised: Is it possible to introduce an operator
$H_i$, $i=0,1,$ independent of the spectral parameter $z$, such
that the equality \hbox{$H_i\psi^{(i)}=z\psi^{(i)}$},
$\psi^{(i)}\in{\cal H}_i$,
implies \hbox{$M_i(z)\psi^{(i)}=0$}\,?
Obviously, having found such an operator
one would reduce the spectral problem for
the transfer-function $M_i(z)$ to the standard spectral problem
for the operator $H_i$ and, thus, the
completeness and basis properties for the eigenvectors of $M_i$
could be studied in terms of the operator $H_i$ referring to
well known facts from operator theory.
A rigorous answer to the above question was found
in~\cite{MotSPbWorkshop}, \cite{MotRem}
in the case where the spectra $\sigma(A_0)$ and $\sigma(A_1)$
of the entries $A_0$ and $A_1$ are separated from each other,
\begin{equation}
\label{SeparCondSigma}
\dist\{\sigma(A_0),\sigma(A_1)\}>0.
\end{equation}
To this end an operator-valued function $V_i(Y_i)$ on the
space of linear operators in $\cH_i$ was constructed
in~\cite{MotSPbWorkshop}, \cite{MotRem} such that
$V_i(Y_i)\psi^{(i)}=V_i(z)\psi^{(i)}$
for any eigenvector $\psi^{(i)}$ corresponding
to an eigenvalue $z$ of the operator $Y_i$.
The desired operator $H_i$ was
searched for as a solution of the operator equation
\begin{equation}
\label{MainOld}
   H_i=A_i+V_i(H_i), \qquad i=0,1.
\end{equation}
Notice that an equation of the form~(\ref{MainOld})
first appeared explicitly in the paper~\cite{BraunMA}
by {\sc M.\,A.\,Braun}.

The solvability of the equation~(\ref{MainOld})
was announced in~\cite{MotSPbWorkshop}
and proved in~\cite{MotRem} under the assumption
$
\|B_{ij}\|_2<\frac{1}{2}\dist\{\sigma(A_0),\sigma(A_1)\}
$
where $\|B_{ij}\|_2$ stands for the Hilbert-Schmidt norm of the
couplings $B_{ij}$. It was found in~\cite{MotRem} that the
problem of constructing the operators $H_i$ is closely related
to the problem of searching for the invariant subspaces of the
matrix $\bH$ which admit representation in the form of graphs
for some bounded $Q_{ji}:\cH_i\to\cH_j$.  The point is that
under the conditions of~\cite{MotSPbWorkshop}, \cite{MotRem}  the
solutions $H_i$, $i=0,1,$ of Eqs.~(\ref{MainOld}) read
$H_i=A_i+B_{ij}Q_{ji}$ while $Q_{ji}$ are contractions.  The
operators $Q_{ji}$ determine a similarity transform reducing the
matrix $\bH$ to the block-diagonal form $\bH'=\diag\{H_0,H_1\}$
(see~\cite{MotRem}).

The idea of the block diagonalization of the $2\times2$ operator
matrices in terms of the invariant subspaces allowing a graph
representation goes back yet to the paper~\cite{Okubo} by {\sc
S.\,Okubo} (regarding applications of Okubo's approach to
particle physics see, e.\,g.,
Refs.~\cite{Okubo}\,--\,\cite{KorchinShebeko}).  In a mathematically
rigorous way this idea was applied by {\sc V.\,A.\,Malyshev} and
{\sc R.\,A.\,Minlos}~\cite{MalyshevMinlos} to a class of
selfadjoint operators in statistical physics. The techniques of
Ref.~\cite{MalyshevMinlos} are restricted to the case where the
norms of the entries $B_{ij}$ are sufficiently small and the
separation condition~(\ref{SeparCondSigma}) holds, too.  {\sc
V.\,M.\,Adamyan} and {\sc H.\,Langer}~\cite{AdL} proved the
existence of invariant subspaces allowing a graph representation
for arbitrary bounded entries $B_{ij}$ however assuming, instead
of the condition~(\ref{SeparCondSigma}), the essentially
different assumption that the spectrum of one of the entries
$A_i$, $i=0,1$ is situated strictly below the spectrum of the
other one, say $\Max\sigma(A_1)<\Min\sigma(A_0)\,.$ Recently,
the result of~\cite{AdL} was extended by {\sc V.\,M.\,Adamyan,
H.\,Langer, R.\,Mennicken} and {\sc J.\,Saurer}~\cite{AdLMSr} to
the case where $\Max\sigma(A_1)\leq\Min\sigma(A_0)\,$ and where
the couplings $B_{ij}$ were allowed to be unbounded operators
such that, for \mbox{$\alpha_0<\Min\sigma(A_0)$}, the product
\mbox{$(A_0-\alpha_0)^{-1/2}B_{01}$} makes sense as a bounded
operator. The mentioned conditions were then somewhat weakened
by {\sc R.\,Mennicken} and {\sc A.\,A.\,Shkalikov}~\cite{MenShk}
in the case of a bounded entry $A_1$ and the same type of
entries $B_{ij}$ as in~\cite{AdLMSr}.  Instead of the explicit
conditions on the spectra of $A_i$, the paper~\cite{MenShk} uses
a condition on the spectrum of the transfer function $M_1(z)$
itself.  One can check that the spectral component $H_i$ of the
matrix $\bH$ constructed in~\cite{AdL}\,--\,\cite{MenShk}
satisfies the equation~(\ref{MainOld}), at least
in the case where for $j\neq i$ the entry $A_j$ is bounded.

In the present work we study the equation~(\ref{MainOld}) in a
case which is totally different from the spectral situations
considered in~\cite{MalyshevMinlos}, \cite{MotSPbWorkshop},
\cite{MotRem}, \cite{AdL}\,--\,\cite{MenShk}.  From the
beginning, we suppose that
\mbox{$\sigma(A_0)\cap\sigma(A_1)\neq\emptyset$.} Moreover, we
are especially interested just in the case where the spectrum
of, say $A_1$, is partly or totally embedded into the continuous
spectrum of $A_0$. We work under the assumption that the
coupling operators $B_{ij}$ are such that the transfer function
$M_1(z)$ admits analytic continuation, as an operator-valued
function, under the cuts along the branches of the absolutely
continuous spectrum $\sigma_{ac}(A_0)$ of the entry $A_0$.  In
Sect.~\ref{Transfer_functions} we describe the conditions on
$B_{ij}$ making such a continuation of $M_1(z)$ possible.

The problem considered is closely related to the resonances
generated by the matrix $\bH$.  Regarding a definition of the
resonance and history of the subject see, e.\,g.,
the book~\cite{ReedSimonIII}. A recent survey of the literature on
resonances can be found in~\cite{MotMathNach}. Throughout the
paper we treat resonances as the discrete spectrum of the
transfer function $M_1(z)$ situated in the unphysical sheets of
its Riemann surface.  We assume that the absolutely continuous
spectrum of the entry $A_0$ consists of $\sm$
($1\leq\sm<\infty$) distinct intervals.  As  a result in Sect.
\ref{SmainEq} we get $2^\sm$ variants of the function $V_1(Y)$
and, consequently, $2^\sm$ different variants of the
equation~(\ref{MainOld}) which read now as Eq.~(\ref{MainEq}).
The solutions of~(\ref{MainEq}) represent non-selfadjoint
operators the spectrum of which includes the resonances in
unphysical sheets neighboring the physical one.

In Sect.~\ref{SecFactor} we first prove the factorization
theorem for the transfer function $M_1(z)$.  It follows from
this theorem that there exist certain domains surrounding the
set $\sigma(A_1)$ and lying partly in the unphysical sheet(s)
where the spectrum of $M_1$ is represented only by the spectrum
of the respective solutions of the basic
equation~(\ref{MainEq}). Since the root vectors of these
solutions are also root vectors for $M_1$, this fact allows us
to talk further, in Sect.~\ref{RealEigen}
and~\ref{BasisnessGeneral}, about completeness and basis properties
of the root vectors of the transfer function $M_1$ corresponding
to its spectrum in the above domains, including the resonance
spectrum. To prove these properties we rely mainly on the
respective statements from the books by {\sc I.\,C.\,Gohberg}
and {\sc M.\,G.\,Krein}~\cite{GK} and by {\sc
T.\,Kato}~\cite{Kato}.

\section{Analytic continuation of the transfer function}
\label{Transfer_functions}

The transfer function $M_i(z)$, $i=0,1$, considered on the
resolvent set $\rho(A_j)$ of the entry $A_j$, $j\neq i$,
represents a particular case of a holomorphic operator-valued
function. In the present work we use the standard definition of
holomorphy of an operator-valued function with respect to the
operator norm topology (see, e.\,g.,~\cite{AdLMSr}).
One can extend to operator-valued functions the usual
definitions of the spectrum and its components.
Each transfer function $M_i(z)$, $i=0,1$,
is holomorphic at least in the resolvent set $\rho(A_j)$
of the entry $A_j$, $j\neq i$. Since the inverse transfer
functions $M_i^{-1}(z)$ coincide with the respective
block components $\bR_{ii}(z)$ of the resolvent
$\bR(z)=({\bf H}-z)^{-1}$,
they are both holomorphic at least in the set $\rho(\bH)$.

Let  $E_j$ be the spectral measure for the entry
$A_j$, $A_j=\int_{\sigma(A_j)}\lambda\,dE_j(\lambda)$\,, $j=0,1,$
\mbox{$\sigma(A_j)\subset\R$}. Then the functions $V_i(z)$
can be written
$$
    V_i(z)=B_{ij}\int_{\sigma(A_j)}dE_j(\mu)\frac{1}{z-\mu}B_{ji}.
$$
Thus, it is convenient to introduce the quantities
$\cV_j(B)=\Sup_{\left\{\delta_k\right\}}
\sum_k \|B_{ij}E_j(\delta_k)B_{ji}\|,$
with $\left\{\delta_k\right\}$ being a finite or countable
complete system of Borel subsets of $\sigma(A_j)$ such
that $\delta_k\cap\delta_l=\emptyset$, if $k\neq l$ and
$\mathop{\bigcup}_k\delta_k=\sigma(A_j)$.  The number
$\cV_j(B)$ is called the variation of the operators $B_{ij}$
with respect to the spectral measure $E_j$.
Along with the ``total'' variation
$\cV_j(B)$ we shall use the ``truncated'' variations
$
    \reduction{\cV_j(B)}{\Delta}=\Sup_{\left\{\delta_k\right\}}
 \sum_k \|B_{ij}E_j(\delta_k\cap\Delta)B_{ji}\|
$
where $\Delta$ is a certain Borel subset of $\sigma(A_j)$;
$\reduction{\cV_j(B)}{\Delta}\leq \cV_j(B)$.

The value of $\cV_j(B)$ satisfies the estimates
$
\|B_{ij}\|^2\leq \cV_j(B)\leq\|B_{ij}\|_2^2.
$
The equality $\cV_j(B)=\|B_{ij}\|^2$ is attained in the case
where $A_j$ is a multiple of the identity operator. The equality
$\cV_j(B)=\|B_{ij}\|_2^2$ holds if $A_j$ possesses only a pure
discrete spectrum which is at the same time simple.

We assume that the spectrum of the operator $A_1$ can intersect
only the continuous spectrum of the operator $A_0$ and this
intersection is only realized on (every of) the pairwise
nonintersecting open intervals
$\Delta_k^0=(\mu_k^{(1)},\mu_k^{(2)})\subset\sigma_c(A_0)$,
$\mu_k^{(1)}<\mu_k^{(2)}$,
\mbox{$k=1,2,\ldots,\sm$,}
\mbox{$\sm<\infty$,}
and \mbox{$-\infty\leq\mu_1^{(1)}$,}
\mbox{$\mu_\sm^{(2)}\leq+\infty$.}
Therefore, we assume that
$\Delta_k^0\cap\sigma(A_1)\neq\emptyset$ for all
$k=1,2,\ldots,\sm$ and $\sigma(A_1)\cap\sigma'(A_0)=\emptyset$
where $\sigma'(A_0)=\sigma(A_0)\setminus
\mathop{\bigcup}_{k=1}^\sm\Delta_k^0$.

We shall suppose that the product
$
  K_B(\mu)\mathop{\mbox{\large$=$}}\limits^{\rm def}
  B_{10}E^0(\mu)B_{01}
$
where $E^0(\mu)$ stands for the spectral function of $A_0$,
$E^0(\mu)=E_0\biggl((-\infty,\mu)\biggr)$,
is differentiable in $\mu\in\Delta_k^0$,
$k=1,2,\ldots,\sm,$ in the operator norm topology.
The derivative $K'_B(\mu)$
is non-negative,
$
 K'_B(\mu)\geq0,
$
since $K_B(\mu)$ is a non-decreasing function.
Differentiability of $K_B(\mu)$ means that the continuous
spectrum of the entry $A_0$ includes, in each $\Delta_k^0$,
$k=1,2,\ldots,\sm,$ a branch of the absolutely continuous
spectrum $\sigma_{ac}(A_0)$.
Obviously,
$
   \reduction{\cV_0(B)}{\Delta_k^0}=\int_{\Delta_k^0}d\mu\,\|K'_B(\mu)\|.
$

Further, we suppose that the function $K'_B(\mu)$ is continuous
within the closed intervals $\overline{\Delta_k^0}$ and,
moreover, that it admits analytic continuation from each of
these intervals to a simply connected domain situated, say, in
$\C^-$.
For the interval $\Delta_k^0$, let this domain be
called $D_k^-$. We assume
that the boundary of each domain $D_k^-$ includes the
entire spectral interval $\Delta_k^0$ and
the domains $D_k^-$ and $D_j^-$ for different $k$ and $j$
do not intersect each other.
Since $K'_B(\mu)$ represents a self-adjoint operator for
$\mu\in\Delta_k^0$ and $\Delta_k^0\subset\R$, the function
$K'_B(\mu)$ also automatically admits analytic continuation from
$\Delta_k^0$ into the domain $D_k^+$, symmetric to $D_k^-$ with
respect to the real axis, $D_k^+=\{z:\,\overline{z}\in D_k^-\}$.
For the continuation into $D_k^+$ we retain the same notation
$K'_B(\mu)$. The selfadjoitness of $K'_B(\mu)$ for
$\mu\in\Delta_k^{0}$ implies
$[K'_B(\mu)]^*=K'_B(\bar{\mu}),$ $\mu\in D_k^\pm\,.$
Also, we shall always suppose the $K'_B(\mu)$
satisfies the H\"older condition at the (finite) end points
$\mu_k^{(1)}$, $\mu_k^{(2)}$ of the spectral intervals
$\Delta_k^0$,
$\|K'_B(\mu)-K'_B(\mu_k^{(i)})\|\leq C|\mu-\mu_k^{(i)}|^\gamma,$
$i=1,2,$ $\mu\in D_k^\pm,$ with some positive $C$ and $\gamma$.

Let $l=(l_1,l_2,\ldots,l_\sm)$ be a multi-index having the
components $l_k=+1$ or $l_k=-1$, $k=1,2,\ldots,\sm$. In what
follows we consider the domains
$D_l=\mathop{\bigcup}_{k=1}^\sm D_k^{l_k}$.
Let $\Gamma_k^{l_k}$ be a rectifiable Jordan curve in
$D_k^{l_k}$ resulting from continuous deformation of the
interval $\Delta^0_k$, the (finite) end points of this interval being
fixed. With the exception of the end
points, the closure $\overline{\Gamma}_k^{l_k}$
of the contour $\Gamma_k^{l_k}$ should have no other common points
with the set $\sigma_c(A_0)$. By $\Gamma_l$,
$l=(l_1,l_2,\ldots,l_\sm)$, we shall denote the union of the
contours $\Gamma_k^{l_k}$,
$
\Gamma_l=\mathop{\bigcup}_{k=1}^\sm\Gamma_k^{l_k}.
$
Also, we extend the definition of the variation $\cV_0(B)$
to the set $\sigma'(A_0)\bigcup\Gamma_l$ by introducing
the modified variation
\begin{equation}
\label{NBnorm}
\cV_0(B,\Gamma_l)=\reduction{\cV_0(B)}{\sigma'(A_0)}+
\displaystyle\int_{\Gamma_l}|d\mu|\,\|K'_B(\mu)\|
\end{equation}
with $|d\mu|$ Lebesgue measure on $\Gamma_l$.  We suppose that
the operators $B_{ij}$ are such that there exists a contour
(contours) $\Gamma_l$ where the value $\cV_0(B,\Gamma_l)$ is
finite, $\cV_0(B,\Gamma_l)<\infty,$ including also the case of
the unbounded set
$\mathop{\bigcup}_{k=1}^\sm\Delta^0_k$.  The contours
$\Gamma_l$ satisfying the condition~$\cV_0(B,\Gamma_l)<\infty$
are said to be {\em $K_B$-bounded} contours.

\begin{lemma}\label{M1-Continuation}
The analytic continuation of the transfer function $M_1(z)$,
$z\in\C\setminus{\sigma(A_0)}$,
through the spectral intervals $\Delta^0_k$ into the subdomain
$D(\Gamma_l)\subset D_l$ bounded by the set
$\mathop{\bigcup}_{k=1}^{\sm}\Delta^0_k$ and
a $K_B$-bounded contour $\Gamma_l$ is given by
\begin{equation}
\label{Mcmpl}
    M_1(z,\Gamma_l)=A_1-z+V_1(z,\Gamma_l)
\end{equation}
where
\begin{eqnarray}
\label{MGamma}
V_1(z,\Gamma_l) &=&\int_{\sigma'(A_0)\cup\Gamma_l}
K_B(d\mu)\,\,\frac{1}{z-\mu}   \\
\nonumber
 &\mathop{\mbox{\large$=$}}\limits^{\rm def} &
\int_{\sigma'(A_0)}B_{10}E_0(d\mu)B_{01}\,\frac{1}{z-\mu} +
\int_{\Gamma_l}
d\mu\,K'_B(\mu)\,\frac{1}{z-\mu}.
\end{eqnarray}
For $z\in D_k^{l_k}\cap D(\Gamma_l)$ the function
$M_1(z,\Gamma_l)$ may be written as
\begin{equation}
\label{M1Gresidue}
  M_1(z,\Gamma_l)=M_1(z)+2\pi\ri\,l_k K'_B(z).
\end{equation}
\end{lemma}

\noindent P~r~o~o~f~.~~The proof is reduced to the observation
that the function $M_1(z,\Gamma_l)$ is holomorphic
for $z\in\C\setminus[\sigma'(A_0)\cup\Gamma_l]$
and coincides with $M_1(z)$ for
$z\in\C\setminus[\sigma'(A_0)\cup\overline{D(\Gamma_l)}]$.
Eq.~(\ref{M1Gresidue}) is obtained from~(\ref{MGamma})
using the Residue Theorem. %
{\nopagebreak\mbox{\phantom{MMMM}}\hfill $\Box$\par\addvspace{0.25cm}}

The formula~(\ref{M1Gresidue}) shows that in general the
transfer function $M_1$ has a multi-sheeted Riemann surface.
The sheet of the complex plane where the transfer function
$M_1(z)$ is considered together with the resolvent $\bR(z)$
initially is said to be the {\em physical sheet}.  The remaining
sheets of the Riemann surface of $M_1(z)$ are said to be {\em
unphysical sheets.}
In the present work we only deal with the unphysical sheets
neighboring the physical one, i.\,e., with the sheets connected
through the intervals $\Delta_k^0$ for some
$k\in\{1,2,\ldots,\sm\}$ immediately to the physical sheet.

\section{The basic equation and its solutions}
\label{SmainEq}

If an operator-valued function
$T:\,\sigma'(A_0)\cup\Gamma\rightarrow\bB(\cH_1,\cH_1)$
is continuous and bounded on a $K_B$-bounded contour
$\Gamma$,
$\|T\|_{\infty,\Gamma}=
\Sup\limits_{\mu\in\sigma'(A_0)\cup\Gamma}\|T(\mu)\|<\infty,
$
and satisfies the Lipschitz condition on $\sigma'(A_0)$
then the integral
\begin{equation}
\label{Xi}
\int_{\sigma'(A_0)\cup\Gamma}K_B(d\mu)\,T(\mu)
\mathop{\mbox{\large$=$}}\limits^{\rm def}
\int_{\sigma'(A_0)}B_{10}E_0(d\mu)B_{01}T(\mu)+
\int_{\Gamma}
d\mu\,K'_B(\mu)\,T(\mu)
\end{equation}
exists in the sense of the operator norm topology
(see  Lemma~7.2 in~\cite{AdLMSr}) and
\begin{equation}
\label{Lem1}
\biggl\|\,\,\int_{\sigma'(A_0)\cup\Gamma}K_B(d\mu)T(\mu)\biggr\|
\leq \cV_0(B,\Gamma)\,\|T\|_{\infty,\Gamma}.
\end{equation}
In particular, if $T(z)$ is the resolvent of an operator
$Y$, \mbox{$T(z)=(Y-z I_1)^{-1}$,} the spectrum of which
is separated from
$\sigma'(A_0)\cup\Gamma$, then one can define the operator
$$
V_1(Y,\Gamma)=
\int_{\sigma'(A_0)\cup\Gamma}
K_B(d\mu)(Y-\mu)^{-1}.
$$
This operator is bounded, $V_1(Y,\Gamma)\in\bB(\cH_1,\cH_1)$,
and, because of~(\ref{Lem1}), its norm admits the estimate
\begin{equation}
\label{V1Yest}
\|V_1(Y,\Gamma)\|\leq \cV_0(B,\Gamma_l)\,
\Sup_{\mu\in\sigma'(A_0)\cup\Gamma}\|(Y-\mu)^{-1}\|.
\end{equation}

In what follows we consider the equation
\begin{equation}
\label{MainEq}
Y=A_1+V_1(Y,\Gamma).
\end{equation}
This equation possesses the following
characteristic property: If an operator $H_1$ is a
solution of~(\ref{MainEq}) and $u_1$ is an eigenvector of $H_1$,
$H_1 u_1=zu_1$, then
automatically
$zu_1=A_1 u_1+V_1(H_1,\Gamma)u_1=A_1 u_1+V_1(z,\Gamma)u_1.$
This implies that any eigenvalue $z$ of such an operator $H_1$
is automatically an eigenvalue for the continued transfer function
$M_1(z,\Gamma_l)$ and $u_1$ is its eigenvector.  Thus, having
found the solution(s) of the equation~(\ref{MainEq}) one obtains an
effective means of studying the spectral properties of the transfer
function $M_1(z,\Gamma)$ itself.
Often it turns out to be convenient to rewrite Eq.~(\ref{MainEq}) in
the form
\begin{equation}
\label{MainEqC}
X=V_1(A_1+X,\Gamma)\,
\end{equation}
where $X=Y-A_1$.
\begin{theorem}\label{Solvability}
Let a contour $\Gamma$ be $K_B$-bounded and
\begin{equation}
\label{Best}
\cV_0(B,\Gamma)< \displaystyle\frac{1}{4}\,d_0^2(\Gamma)\,
\end{equation}
where
$d_0(\Gamma)=\dist\{\sigma(A_1),\sigma'(A_0)\cup\Gamma\}.$
Then Eq.~{\rm(\ref{MainEqC})} is uniquely solvable
in any ball
${\cal S}_1(r)\subset\bB(\cH_1,\cH_1)$ including operators
$X:\cH_1\rightarrow\cH_1$ the norms of which are bounded as
$\|X\|\leq r$ with $r$ such that
\begin{equation}
\label{Br}
  r_{\rm min}(\Gamma)\leq r < r_{\rm max}(\Gamma).
\end{equation}
Here
\begin{equation}
\label{rminmax}
r_{\rm min}(\Gamma)=d_0(\Gamma)/2- \sqrt{d_0^2(\Gamma)/4
-\cV_0(B,\Gamma)}, \qquad
r_{\rm max}(\Gamma)=d_0(\Gamma)-\sqrt{\cV_0(B,\Gamma)}.
\end{equation}
The solution $X$ of Eq.~{\rm(\ref{MainEqC})} is the same for any
$r$ satisfying~{\rm(\ref{Br})} and in fact it belongs to the smallest
ball ${\cal S}_1(r_{\rm min})$, $\|X\|\leq r_{\rm min}(\Gamma)$.
\end{theorem}
\noindent P~r~o~o~f~.~
Let $F(X)=V_1(A_1+X,\Gamma)$ with $X{\in}{\cal S}_1(r)$.
To begin with we search for
a condition under which the function $F$ is a contracting mapping of
the ball ${\cal S}_1(r)$ into itself.  Since, in view of~(\ref{Br}) and
(\ref{rminmax}) the condition $0<r<d_0$, $d_0=d_0(\Gamma)$
automatically holds, the spectrum of the operator $A_1+X$ does
not intersect the set $\sigma'(A_0)\cup\Gamma$.
This means that for all
$\mu\in\sigma'(A_0)\cup\Gamma$ the resolvent
$(A_1+X-\mu{I_1})^{-1}$ exists as a bounded operator in
$\cH_1$. It follows from the estimate~(\ref{V1Yest}) that
$$
\|F(X)\|\leq \cV_0(B,\Gamma)
\Sup_{\mu\in\sigma'(A_0)\cup\Gamma}\|(A_1+X-\mu)^{-1}\|
\leq\cV_0(B,\Gamma)\,\displaystyle\frac{1}{d_0-r}
$$
while
$$
\|F(X)-F(Y)\|\leq
\cV_0(B,\Gamma)\,\displaystyle\frac{1}{(d_0-r)^2}\|Y-X\|.
$$
It immediately follows from these
that the ball ${\cal S}_1(r)$ is necessarily mapped
by the function $F$ into itself and this mapping is a
contraction if the radius $r$ and the value $\cV_0(B,\Gamma)$
are such that
\begin{equation}
\label{Est1}
\cV_0(B,\Gamma)\leq r(d_0-r),\qquad\cV_0(B,\Gamma)<(d_0-r)^2\,.
\end{equation}
Under the condition~(\ref{Best}) the inequalities~(\ref{Est1})
considered together are just equivalent
to the condition~(\ref{Br}).
Thus if this condition is valid then
Eq.~(\ref{MainEqC}) has a solution in any ball ${\cal S}_1(r)$
with $r$ satisfying~(\ref{Br}) and this solution is unique.
This means that the solution is the same
for all the radii satisfying~(\ref{Br}). Moreover, it belongs to the
ball ${\cal S}_1(r_{\rm min})$ with the radius $r_{\rm min}$
given by~(\ref{rminmax}).%
{\nopagebreak\mbox{\phantom{MMMM}}\hfill $\Box$\par\addvspace{0.2cm}}
\begin{theorem}\label{Hunique}
Let the conditions of {\rm Theorem~\ref{Solvability}} be valid
for a $K_B$-bounded contour $\Gamma\subset D_l$
and let $X$ be the solution of Eq.{~\rm(\ref{MainEqC})}
referred to there. Then $X$ coincides with the analogous solution
$\tilde{X}$ for any other $K_B$-bounded
contour $\tilde{\Gamma}\subset D_l$ satisfying the estimate
$\cV_0(B,\tilde{\Gamma})<\tilde{d}_0^2/4$ with
$0 < \tilde{d}_0=\dist\{\sigma(A_1),
\sigma'(A_0)\cup\tilde{\Gamma}\}\leq d_0(\Gamma)$.
Moreover, this solution satisfies the inequality
$\|X\|\leq r_0(B)$ wherem
\hbox{$
%
  r_0(B)=\Inf\limits_{\Gamma_l:\,\omega(B,\Gamma_l)>0}
  r_{\rm min}(\Gamma_l)\,
$}
with $r_{\rm min}(\Gamma_l)$ given by~{\rm(\ref{rminmax})} and
$\omega(B,\Gamma_l)=d_0^2(\Gamma_l)-4\cV_0(B,\Gamma_l).$ The value of
$r_0(B)$ does not depend on $l$.
\end{theorem}
\noindent P~r~o~o~f of this theorem is reduced to an appropriate
continuous deformation of the integration paths.
An only essential point is a checking
of independence of the radius
$r_0(B)$ of the multi-index $l$.  To be sure in such an independence
we consider an arbitrary
$K_B$$-$bounded contour
$\Gamma_l\subset D_l$,
$\Gamma_l=\mathop{\bigcup}_{k=1}^\sm\Gamma_k^{l_k}$.
Denote by $\Gamma_{l'}$ a contour resulting from $\Gamma_l$
after replacing a part of $\Gamma_k^{l_k}$ with
curves
\mbox{$\Gamma_k^{(-l_k)}=\{\mu:\,\overline{\mu}\in\Gamma_k^{l_k}\}$},
symmetric to $\Gamma_k^{l_k}$ with respect to the real axis.
Obviously, such replacements generate,
additionally to
$\Gamma_l$, \mbox{$2^\sm-1$} different contours $\Gamma_{l'}$
where \mbox{$l'=(l'_1,l'_2,\ldots,l'_\sm)$} with $l'_k=\pm l_k$,
$k=1,2,\ldots,\sm.$ For every of such contours
the value of
$\cV_0(B,\Gamma_{l'})$ is the same, namely
\hbox{$\cV_0(B,\Gamma_{l'})=\cV_0(B,\Gamma_l)$,}
since the replacement of $\Gamma_l$ with $\Gamma_{l'}$ does not change
\hbox{$\displaystyle\int_{\Gamma_l}|d{\mu}|\,\|K'_B({\mu})\|$}.
But this just means that
$r_0(B)$ does not depend on $l$.%
{\nopagebreak\mbox{\phantom{MMMM}}\hfill
$\Box$\par\addvspace{0.2cm}}

So, for a given holomorphy domain $D_l$ the solutions $X$ and
$H_1$, $H_1=A_1+X,$ do not depend on the $K_B$-bounded contours
$\Gamma_l\subset D_l$ satisfying the condition~(\ref{Best}).
But when the index $l$ changes,  $X$ and $H_1$ can also change.
For this reason we shall supply them in the following, when it
is necessary, with the index $l$ writing, respectively,
$X^{(l)}$ and $H_1^{(l)}$, $H_1^{(l)}=A_1+X^{(l)}$.  Therefore,
Theorem~\ref{Solvability} guarantees us, in general, the
existence of the $2^\sm$ solutions $X^{(l)}$ to the basic
equation~(\ref{MainEq}) and, hence, the $2^\sm$ respective
solutions $H_1^{(l)}$ to the basic equation~(\ref{MainEqC}).
Surely, Eqs.~(\ref{MainEq}) and~(\ref{MainEqC}) are non-linear
equations and, outside the balls $\|X\|<r_{\rm max}(\Gamma)$,
they may, in principle, have other solutions, different from the
$X^{(l)}$ or $H_1^{(l)}$ the existence of which is guaranteed by
Theorem~\ref{Solvability}. In the following we shall deal only
with the solutions $X^{(l)}$ or $H_1^{(l)}$.
\section{Factorization of the transfer function}
\label{SecFactor}
\begin{theorem}\label{factorization}
Let $\Gamma_l$ be a $K_B$-bounded contour satisfying the
condition~{\rm(\ref{Best})} and $H_1^{(l)}=A_1+X^{(l)}$ with
$X^{(l)}$ the above solution of the basic
equation~{\rm(\ref{MainEq})}, $\|X^{(l)}\|\leq r_0(B)$. Then,
for $z\in\C\setminus(\sigma'(A_0)\cup\Gamma_l)$, the
transfer function $M_1(z,\Gamma_l)$ admits the factorization
\begin{equation}
\label{Mfactor}
    M_1(z,\Gamma_l)=W_1(z,\Gamma_l)\,(H_1^{(l)}-z)\,
\end{equation}
where $W_1(z,\Gamma_l)$ is a bounded operator in $\cH_1$,
\begin{equation}
\label{Mtild}
 W_1(z,\Gamma_l)=I_1-\int_{\sigma'(A_0)\cup\Gamma_l}
 K_B(d\mu)\frac{1}{\mu-z}(H_1^{(l)}-\mu)^{-1}\,.
\end{equation}
Here, $I_1$ stands for the identity operator in $\cH_1$.

For $\dist\{z,\sigma(A_1)\}\leq{d_0(\Gamma_l)/2}$ the
operator $W_1(z,\Gamma_l)$ is boundedly invertible and
\begin{equation}
\label{Mtest}
  \left\|[W_1(z,\Gamma_l)]^{-1}\right\|
 \leq \left(1-\frac{\cV_0(B,\Gamma_l)}
  {d_0^2(\Gamma_l)/4}\right)^{-1} <\infty.
\end{equation}
\end{theorem}

Note that the above statement recalls the known
factorization theorem by {\sc A.\,I.\,Vi\-ro\-zub} and {\sc
V.\,I.\,Matsaev}~\cite{VirozubMatsaev} being valid for a class
of selfadjoint operator-valued functions
(see also~\cite{MarkusMatsaev}).  However, in the case we
deal with in the present paper, the function $M_1(z,\Gamma_l)$
it is not even a selfadjoint operator-valued function in the
sense of~\cite{VirozubMatsaev}.

\bigskip

\noindent P~r~o~o~f~.  For
$z\in\C\setminus(\sigma'(A_0)\cup\Gamma_l)$,  the boundeness of
the operator $W_1(z,\Gamma_l)$ given
by~(\ref{Mtild}) is evident. To prove the
factorization~(\ref{Mfactor}) we note that
for any $z\not\in\sigma'(A_0)\cup\Gamma_l$
\begin{equation}
\label{Th2Start}
 W_1(z,\Gamma_l)\,(H_1^{(l)}-z)=H_1^{(l)}-z
-\int_{\sigma'(A_0)\cup\Gamma_l}
 K_B(d\mu)\frac{(H_1^{(l)}-\mu)^{-1}(H_1^{(l)}-z)}{\mu-z}.
\end{equation}
Since
$(H_1^{(l)}-\mu)^{-1}(H_1^{(l)}-z)=I_1+(\mu-z)(H_1^{(l)}-\mu)^{-1}$
one finds
$$
 \int_{\sigma'(A_0)\cup\Gamma_l}
 K_B(d\mu)\frac{(H_1^{(l)}-\mu)^{-1}(H_1^{(l)}-z)}{\mu-z} =
\int_{\sigma'(A_0)\cup\Gamma_l} \frac{K_B(d\mu)}{\mu-z}
+\int_{\sigma'(A_0)\cup\Gamma_l} K_B(d\mu)(H_1^{(l)}-\mu)^{-1}.
$$
But according to~(\ref{Mcmpl})
$
\int_{\sigma'(A_0)\cup\Gamma_l}
K_B(d\mu)(\mu-z)^{-1}=A_1-z-M_1(z,\Gamma_l),
$
while according to~(\ref{MainEq})
$
\int_{\sigma'(A_0)\cup\Gamma_l}
 K_B(d\mu)(H_1^{(l)}-\mu)^{-1}=H_1^{(l)}-A_1.
$
Making use of these expressions one immediately obtains
Eq.~(\ref{Mfactor}).

Further, we prove that the factor $W_1(z,\Gamma_l)$ is a
boundedly invertible operator if the condition
$\dist\{z,\sigma(A_1)\}\leq{d_0(\Gamma_l)/2}$ is valid. Indeed, under this
condition
$
|\mu-z|\geq\dist\{z,\sigma'(A_0)\cup\Gamma_l\}\geq
d_0(\Gamma_l)/2
$
since
$\dist\{\sigma(A_1),\sigma'(A_0)\cup\Gamma_l\}=d_0(\Gamma_l)$.
On the other hand
$H_1^{(l)}=A_1+X^{(l)}$ and $\|X^{(l)}\|<d_0(\Gamma_l)/2$.
Thus for $\mu\in\sigma'(A_0)\cup\Gamma_l$ we have
$
\|(H_1^{(l)}-\mu)^{-1}\| <2/d_0(\Gamma_l).
$
Consequently
\begin{equation}
\label{forOm}
\biggl\|\, \int_{\sigma'(A_0)\cup\Gamma_l}
 K_B(d\mu)\frac{1}{\mu-z}(H_1^{(l)}-\mu)^{-1}\biggr\|
< \frac{\cV_0(B,\Gamma_l)}{(d_0(\Gamma_l)/2)^2}<1
\end{equation}
and, thus, the estimate~(\ref{Mtest}) is true.
{\nopagebreak\mbox{\phantom{MMMM}}\hfill $\Box$\par\addvspace{0.2cm}}

It is easy to write some simple but useful relations between a
part of the operators $H_1^{(l)}$.  In particular, we derive
such relations between $H_1^{(l)}$ and $H_1^{(-l)}$,
$(-l)=(-l_1,-l_2,\ldots,-l_\sm)$ where $l_k$,
$k=1,2,\ldots,\sm,$ stand for the components of the multi-index
$l=(l_1,l_2,\ldots,l_\sm)$. According to our convention,
$\Gamma_{(-l)}$, $\Gamma_{(-l)}\subset D_{(-l)}$, is a contour
which is obtained from the contour $\Gamma_l$ by replacing all
the components $\Gamma_k^{l_k}$ with the conjugate ones
$\Gamma_k^{(-l_k)}$.
\begin{lemma}\label{Adjoint}
Let $\Gamma_l\subset D_l$ be a $K_B$-bounded contour
for which the conditions of {\rm Theorem~\ref{Solvability}}
are valid. Then for any
$z\in\C\setminus\biggl(\sigma'(A_0)\cup\Gamma_l\biggr)$
the following equality holds true:
\begin{equation}
\label{Hadj}
W_1(z,\Gamma_l)\,\biggl(H_1^{(l)}-z\biggr)=
\biggl(H_1^{(-l)*}-z\biggr)\,
[W_1(\overline{z},\Gamma_{(-l)})]^*\,.
\end{equation}
Therefore the spectrum of $H_1^{(-l)*}$
coincides with the spectrum of $H_1^{(l)}$.
\end{lemma}
\begin{theorem}\label{SpHalfVic}
The spectrum $\sigma(H_1^{(l)})$ of the operator
$H_1^{(l)}=A_1+X^{(l)}$ belongs to the closed $r_0(B)$-vicinity
${\cal O}_{r_0}(A_1)$ of the spectrum of $A_1$, ${\cal
O}_{r_0}(A_1)= \{z\in\C:\,\dist\{z,\sigma(A_1)\}\leq r_0(B)\}$.
If a contour $\Gamma_l\subset D_l$ satisfies~{\rm(\ref{Best})},
then the complex spectrum of $H_1^{(l)}$ belongs  to
$D_l\cap{\cal O}_{r_0}(A_1)$ while outside $D_l$ the spectrum
of $H_1^{(l)}$ is pure real.  Moreover, the spectrum
$\sigma(H_1^{(l)})$ coincides with a (subset of the) spectrum of the
transfer function $M_1(z,\Gamma_l)$. More precisely, the spectrum of
$M_1(z,\Gamma_l)$ in ${\cal O}_{d_0/2}(A_1)=\{z:\, z\in\C,
\dist\{z,\sigma(A_1)\}\leq{d_0(\Gamma_l)}/{2}\}$ is represented
only by the spectrum of $H_1^{(l)}$, i.\,e.
$\sigma\biggl(M_1(\cdot,\Gamma_l)\biggr)\cap{\cal
O}_{d_0/2}(A_1)=\sigma(H_1^{(l)})$. In fact such a statement is true
separately for point and continuous spectra.
\end{theorem}
\begin{theorem}\label{HlpHl2p}
Suppose that two different domains $D_{l'}$ and $D_{l''}$
include the same subdomain $D_k^{l_k}$ for some
$k=1,2,\ldots,\sm$, i.\,e., $l'_k=l''_k=l_k$.
Then the spectra of the operators $H_1^{(l')}$
and $H_1^{(l'')}$ in $D_k^{l_k}$ coincide.
\end{theorem}
Let
\hbox{$
\Omega^{(l)}=\displaystyle\int_{\sigma'(A_0)\cup\Gamma_l}
(H_1^{(-l)*}-\mu)^{-1}K_B(d\mu)\,
(H_1^{(l)}-\mu)^{-1}\,,
$}
where as previously
where $\Gamma_l$ stands for a $K_B$-bounded
contour satisfying the condition~(\ref{Best}).
The operator
$\Omega^{(l)}$ does not depend (for a fixed $l$) on the choice
of such a $\Gamma_l$.
At the same time
$\Omega^{(-l)}=\Omega^{(l)*}.$
The norm $\Omega^{(l)}$ satisfies the estimate
\begin{equation}
\label{Omest}
\|\Omega^{(l)}\|<\frac{\cV_0(B,\Gamma_l)}{(d_0(\Gamma_l)/2)^2}<1.
\end{equation}
\begin{theorem}
\label{MHOmega}
The operators $\Omega^{(l)}$ possess the following properties%
{\rm(}cf. {\rm\cite{MenShk}, \cite{VirozubMatsaev}, \cite{MarkusMatsaev}):}
\begin{eqnarray*}
%
-\frac{1}{2\pi\ri}\int_\gamma dz\,[M_1(z,\Gamma_l)]^{-1} &=&
(I_1+\Omega^{(l)})^{-1}\,,\\
%
%
-\frac{1}{2\pi\ri}\int_\gamma dz\,z\,[M_1(z,\Gamma_l)]^{-1} &=&
(I_1+\Omega^{(l)})^{-1}H_1^{(-l)*}=
H_1^{(l)}(I_1+\Omega^{(l)})^{-1}\,,
\end{eqnarray*}
where $\gamma$ stands for an arbitrary rectifiable closed
{\rm(}including the points at infinity if the entry $A_1$ is
unbounded{\rm)} contour going in the positive direction around the
spectrum of $H_1^{(l)}$ inside the set ${\cal
O}_{d_0(\Gamma)/2}(A_1)$.  The integration over $\gamma$ is
understood in the strong sense.
\end{theorem}
\section{Properties of real eigenvalues}
\label{RealEigen}
If $\lambda$ is a real isolated eigenvalue of the operator
$H_1^{(l')}=A_1+X^{(l')}$,
$l'=(l'_1,l'_2,\ldots,l'_\sm)$, then it is such an eigenvalue
also for the remaining $2^{m-1}$
operators $H_1^{(l)}=A_1+X^{(l)}$ where
$l=(l_1,l_2,\ldots,l_\sm)$ with arbitrary $l_k=\pm1$,
$k=1,2,\ldots,\sm.$ The resolvents of every of the $2^\sm$ operators
$H_1^{(l)}$ at $z=\lambda$ have a first order pole.
Simultaneously such an eigenvalue $\lambda$ belongs to
the point spectrum of the total operator $\bH$.

An isolated real eigenvalue $\lambda$ of the operator $H_1^{(l)}$
can not belong to the spectrum $\sigma'(A_0)$ of the
entry $A_0$ lying outside
$\mathop{\bigcup}_{k=1}^\sm\Delta_k^0$.  Indeed,
according to Theorem~\ref{SpHalfVic}, the spectrum of $H_1^{(l)}$
for arbitrary $l$ is situated in the $r_0(B)$-vicinity ${\cal
O}_{r_0}(A_1)$ of the set $\sigma(A_1)$ and in any case
$r_0(B)<\frac{1}{2}\dist\{\sigma'(A_0),\sigma(A_1)\}$  so that
automatically
\hbox{$\sigma'(A_0)\cap\sigma(H_1^{(l)})=\emptyset.$}
Hence, such a $\lambda$ belongs either to the resolvent set
$\rho(A_0)$ of the entry $A_0$ or it is embedded into
the continuous spectrum of $A_0$ in
$\mathop{\bigcup}_{k=1}^\sm\Delta^0_k$.
\begin{lemma}\label{LReal2}
If a vector $\psi^{(1)}\in\cD(A_1)$ is an eigenvector of
$H_1^{(l)}$ corresponding to a real eigenvalue
$\lambda\in\rho(A_0)$
then the vector $\Psi=(\psi^{(0)},\psi^{(1)})\in\cH$ with
\begin{equation}
\label{psi0}
\psi^{(0)}=-R_0(\lambda)B_{01}\psi^{(1)}
\end{equation}
is an eigenvector of $\bH$, $\bH\Psi=\lambda\Psi$. The converse
statement is also true: if $\lambda$, $\lambda\in\rho(A_0)$,
is a real eigenvalue of $H_1^{(l)}$ and
$\bH\Psi=\lambda\Psi$ for some $\Psi=(\psi^{(0)},\psi^{(1)})$
with $\psi^{(0)}\in\cD(A_0)$ and $\psi^{(1)}\in\cD(A_1)$,
then $\psi^{(0)}$ is related to $\psi^{(1)}$ as in~{\rm(\ref{psi0})}
and $H_1^{(l)}\psi^{(1)}=\lambda\psi^{(1)}$.
\end{lemma}

If an eigenvalue $\lambda$ of $H_1^{(l)}$ belongs to
$\Delta_k^0=(\mu_k^{(1)},\mu_k^{(2)})$ for some
$k=1,2,\ldots,\sm$, then
$
|\lambda-\mu_k^{(i)}|\geq\dist\{\mu_k^{(i)},\sigma(A_1)\}-r_0(B),
\quad i=1,2.
$
Recall that according to our assumption the entry $A_0$
has no point spectrum inside $\Delta_k^0$. Since
$\Delta_k^0$ is a part of the continuous spectrum of $A_0$,
the resolvent $R_0(z)=(A_0-z)^{-1}$ for $z=\lambda\pm\ri0$
exists being however an unbounded operator. Nevertheless
a statement analogous to Lemma~\ref{LReal2} is valid in this
case, too.
\begin{lemma}\label{LReal3}
If a vector $\psi^{(1)}\in\cD(A_1)$ is an eigenvector of
$H_1^{(l)}$ corresponding to a real eigenvalue
$\lambda\in\Delta^0_k=(\mu^{(1)}_k,\mu^{(2)}_k)$,
$k=1,2,\ldots,\sm$, \, $H_1^{(l)}\psi^{(1)}=\lambda\psi^{(1)}$,
then either

a{\rm)} $E^0(\mu)B_{01}\psi^{(1)}=0$ for all $\mu\leq\mu_k^{(2)}$

\noindent or

b{\rm)} $E^0(\mu)B_{01}\psi^{(1)}\neq0$ for any $\mu\in\Delta^0_k$,

c{\rm)} the function $\|E^0(\mu)B_{01}\psi^{(1)}\|$ is
differentiable in $\mu$ on $\Delta^0_k$

\noindent and

d{\rm)} $\reduction{\displaystyle\frac{d}{d\mu}
         \|E^0(\mu)B_{01}\psi^{(1)}\|}
         {\mu=\lambda}=0.$

In both cases the vector $\psi^{(0)}$ given by~{\rm(\ref{psi0})}
exists in $\cD(A_0)$ and $\Psi=(\psi^{(0)},\psi^{(1)})$
is an eigenvector of $\bH$, $\bH\Psi=\lambda\Psi.$

The converse statement is also true. Namely, if a
$\Psi=(\psi^{(0)},\psi^{(1)})$ with $\psi^{(0)}\in\cD(A_0)$ and
$\psi^{(1)}\in\cD(A_1)$ is an eigenvector of $\bH$,
$\bH\Psi=\lambda\Psi$, corresponding to an eigenvalue
$\lambda$ of $H_1^{(l)}$, $\lambda\in\Delta^0_k$, then either
the condition {\rm(}a{\rm)} is valid or the conditions
{\rm(}b\,--\,d{\rm)} are valid.  In both cases the
relation~{\rm(\ref{psi0})} is retained meaning, in particular,
that $\psi^{(1)}\neq0$ and $\psi^{(1)}$ is an eigenvector of
$H_1^{(l)}$ corresponding to the eigenvalue $\lambda$.
\end{lemma}
\medskip

Let $\sigma_{pri}(H_1^{(l)})$ be the set of all real isolated
eigenvalues of the operator $H_1^{(l)}$.
As we already established, this set coincides with the part
$\sigma_{pri}(M_1(\cdot,\Gamma_l)$ of the set of the real
isolated eigenvalues of the transfer function $M_1(z,\Gamma_l)$
belonging to ${\cal O}_{d_0/2}(A_1)$ for any  $K_B$-bounded
contour $\Gamma_l$ satisfying the condition~(\ref{Best}).

Since in the remainder of the Section we will consider different
eigenvalues $\lambda\in\sigma_{pri}(H_1^{(l)})$, we will use a
more specific notation, $\psi^{(1)}_{\lambda,j}$,
$j=1,2,\ldots,m_\lambda$, for the respective eigenvectors of the
$H_1^{(l)}$. The notation $m_\lambda$, $m_\lambda\leq\infty,$
stands for the multiplicity of the eigenvalue $\lambda$. Recall
that every $\psi^{(1)}_{\lambda,j}$ is an eigenvector
simultaneously for all the $H_1^{(l)}$ and
$M_1(\lambda\pm\ri0,\Gamma_l)$,
\mbox{$l=(l_1,l_2,\ldots,l_\sm)$} with  $l_k=\pm1$,
$k=1,2,\ldots,\sm$.
In the considered case
the multiplicity $m_\lambda$ is both
the geometric and algebraic multiplicity of $\lambda$.
Respective eigenvectors of the total
matrix $\bH$ will be denoted by $\Psi_{\lambda,j}$,
$\Psi_{\lambda,j}=(\psi^{(0)}_{\lambda,j},\psi^{(1)}_{\lambda,j})$.
It will be supposed that the $\psi^{(1)}_{\lambda,j}$ are chosen
in such a way that the vectors $\Psi_{\lambda,j}$ are
orthonormal, $\lal\Psi_{\lambda,j},\Psi_{\lambda',j'}\ral=
\delta_{\lambda\lambda'}\delta_{jj'}$.  Obviously, the
statements of Lemmas~\ref{LReal2} and~\ref{LReal3} imply that
the eigenvectors $\Psi_{\lambda,j}$,
$\lambda\in\sigma_{pri}(H_1^{(l)})$, $j=1,2,\ldots,m_\lambda,$
form an orthonormal basis in the invariant subspace of the
operator $\bH$ corresponding to the subset
$\sigma_{pri}(H_1^{(l)})$ of the point spectrum $\sigma_p(\bH)$
of $\bH$.
\begin{lemma}\label{PosOmegaReal}
Let $\cH_1^{(pri)}$, $\cH_1^{(pri)}\subset\cH_1$, be the closed
span of the eigenvectors $\psi^{(1)}_{\lambda,j}$< of $H_1^{(l)}$
corresponding to the spectrum $\sigma_{pri}(H_1^{(l)})$,
$
\cH_1^{(pri)}=\overline{{\sf V}\{\psi^{(1)}_{\lambda,j},\,
\lambda\in\sigma_{pri}(H_1^{(l)}),\, j=1,2,\ldots,m_\lambda\}}.
$
For any $l=(l_1,l_2,\ldots,l_\sm)$, $l_k=\pm1$,
$k=1,2,\ldots,\sm$, the operator $\Omega^{(l)}$ is non-negative
on the subspace $\cH_1^{(pri)}$.
\end{lemma}
This statement implies that one can introduce a new inner
product in $\cH_1^{(pri)}$,
{$[u_1,v_1]_{\cH_1^{(pri)}}=\lal(I_1+\Omega^{(l)})u_1,v_1\ral.$}
Then, with the help of a theorem of {\sc N.\,K.\,Bari}
(Theorem~VI.2.1 of~\cite{GK}) one can prove the following
statement.

\begin{theorem}\label{RealRieszBasis}
The system of vectors $\psi^{(1)}_{\lambda,j}$,
$\lambda\in\sigma_{pri}(H_1^{(l)}),$ $j=1,2,\ldots,m_\lambda,$
forms a Riesz basis of the subspace $\cH_1^{(pri)}$.
\end{theorem}
\section{Completeness and basis properties}
\label{BasisnessGeneral}
We restrict ourselves to the case where the entry $A_1$ has pure
discrete spectrum only, i.\,e., the resolvent
$R_1(z)=(A_1-z)^{-1}$ is a compact operator in $\cH_1$ for any
$z\in\rho(A_1)$. In this case the operators $H_1^{(l)}$ have
compact resolvents, too.  This is a consequence of
Theorem~V.3.17 of~\cite{Kato} since the difference
\mbox{$H_1^{(l)}-A_1=X^{(l)}$} is a bounded operator (see
Theorem~\ref{Solvability}). Also, the operators $X^{(l)}$ are
compact.

Denote by $\cH_{1,\lambda}^{(l)}$ the algebraic eigenspace of
$H_1^{(l)}$ corresponding to an eigenvalue $\lambda$.
Let $m_\lambda$ be the algebraic
multiplicity, $m_\lambda=\mathop{\rm dim}\cH_{1,\lambda}^{(l)}$,
$m_\lambda<\infty$, and
$\psi^{(l)}_{\lambda,i},$ $i=1,2,\ldots,m_\lambda,$ be the root
vectors of $H_1^{(l)}$ forming a basis
of the subspace $\cH_{1,\lambda}^{(l)}$.
In the following we will try to give
an answer on the question when the union of such bases in
$\lambda$ forms a basis of the total space $\cH_1$. But, in
any case, we already have  an assertion regarding
completeness of the system
\begin{equation}
\label{RootVecSystem}
\{\psi^{(l)}_{\lambda,i},\,\,\lambda\in\sigma(H_1^{(l)}),\,\,
                              i=1,2,\ldots,m_\lambda\}.
\end{equation}
\begin{theorem}\label{Completeness}
The closure of the linear span of the system~{\rm(\ref{RootVecSystem})}
coincides with $\cH_1$.
\end{theorem}
\noindent This assertion is a particular case
of Theorem~V.10.1 from~\cite{GK}.

We shall consider the case where the intersection
$\biggl(\mathop{{\bigcup}}_{k=1}^{\sm}\Delta^0_k\biggr)
\cap\sigma(A_1)$ includes infinitely many points and the entry
$A_1$ is semibounded from below.
This assumption means that at least the interval
$\Delta_\sm^0$ is infinite,
$\Delta_\sm^0=(\mu_\sm^{(1)},+\infty)$.
The eigenvalues
$\lambda_i^{(A_1)}$, $i=1,2,\ldots\,$, of the operator $A_1$
will be enumerated in increasing order,
$\lambda_1^{(A_1)}<\ldots<\lambda_i^{(A_1)}
<\lambda_{i+1}^{(A_1)}<\ldots\,$\,;
\mbox{$\Lim_{i\to\infty}\lambda_i^{(A_1)}=+\infty$} exists.

Suppose further that there is a number $i_0$ such that
for any $i\geq i_0$ and for some fixed $r>r_0(B)$
\begin{equation}
\label{DifLambda}
  \lambda_i^{(A_1)}-\lambda_{i-1}^{(A_1)}>2r.
\end{equation}
Let $\gamma_0$ be a circle centered at
\mbox{$z=(\lambda_1^{(A_1)}+\lambda_{i_0-1}^{(A_1)})/2$} and
having the radius
\mbox{$(\lambda_{i_0-1}^{(A_1)}-\lambda_1^{(A_1)})/2+r$} while
the $\gamma_i$ for $i\geq i_0$ are the circles with centers
$\lambda_i^{(A_1)}$ and the radius $r$.  Let us introduce the
projections
\mbox{$\sQ_i^{(l)}=-\frac{1}{2\pi\ri}\int_{\gamma_i} dz\,
(H_1^{(l)}-z)^{-1},$} $i=0,i_0,i_0+1,\ldots$\,.  Every
projection $\sQ_i^{(l)}$ represents a sum of the
eigenprojections corresponding to the
eigenvalues $\lambda^{(l)}$ of $H_1^{(l)}$ lying inside
$\gamma_i$ and
\hbox{$\sQ^{(l)}_i\sQ^{(l)}_j=\delta_{ij}\sQ^{(l)}_i.$} The
subspaces \mbox{$\cN_i^{(l)}=\sQ_i^{(l)}\cH_1$} are invariant
under $H_1^{(l)}$; $\mathop{\rm dim}\cN_i^{(l)}$ coincides with a
sum of algebraic multiplicities of the eigenvalues
$\lambda^{(l)}$ lying inside $\gamma_i$.
We introduce also the (orthogonal) projections
\mbox{$\sP_i^{(A_1)}=-\frac{1}{2\pi\ri}\int_{\gamma_i} dz\,
(A_1-z)^{-1},$} $i=0,i_0,i_0+1,\ldots\,\,\,.$
\begin{lemma}\label{NomegaLInd}
Under the condition~{\rm(\ref{DifLambda})} the sequence
$\cN_i^{(l)}$, \, $i=0,i_0,i_0+1,\ldots,$ \,
is $\omega$-linearly independent
and complete in $\cH_1$.
If instead of~{\rm(\ref{DifLambda})} the condition
\begin{equation}
\label{DifLambda4}
  \lambda_i^{(A_1)}-\lambda_{i-1}^{(A_1)}>2r>4r_0(B)\qquad
\forall i\geq i_0,
\end{equation}
is satisfied then
\mbox{$\mathop{\rm dim}\cN_i^{(l)}=
\mathop{\rm dim}\sP_i^{(A_1)}\cH_1$,}
$i=0,i_0,i_0+1,\ldots$\,\,.
\end{lemma}
\begin{theorem}\label{KatoV4-15-16}
Assume $\lambda_{i+1}^{(A_1)}-\lambda_i^{(A_1)}\to\infty$
as $i\to\infty.$ Let $i_0$ be a number starting from which
the inequality
{\rm(\ref{DifLambda4})} holds.
Then the following limit exists
\begin{equation}
\label{sLimQ}
\mathop{s-\rm lim}\limits_{n\to\infty}
\sum_{i=0,i\geq i_0}^n \sQ_i^{(l)}=I_1.
\end{equation}
Additionally, assume that
\begin{equation}
\label{InvSquareConvergent}
\sum_{i=1}^\infty
{(\lambda_{i+1}^{(A_1)}-\lambda_i^{(A_1)})^{-2}}<\infty\,.
\end{equation}
Then~{\rm(\ref{sLimQ})} is true for any renumbering of
$\sQ_i^{(l)}$. Moreover, there exists a constant $C$ such that
$\biggl\|\sum_{i\in\cI}\sQ_i^{(l)}\biggr\|\leq C$ for any finite
set $\cI$ of integers $i=0$, $i\geq i_0$.
\end{theorem}
This theorem represents a slightly extended statement of
Theorems~V.4.15 and~V.4.16 of~\cite{Kato}.
\begin{remark}\label{KatoV4-15-16-Rem1}
Eq.~{\rm(\ref{sLimQ})} implies that
\begin{equation}
\label{sLimP}
\mathop{s-\rm lim}\limits_{n\to\infty}
\sum_{i=0,i\geq i_0}^n\,\,\,
\sum_{\lambda\in\Inter\gamma_i} \sP_\lambda^{(l)}=I_1,
\end{equation}
where $\lambda$ stand for the eigenvalues of the
operator $H_1^{(l)}$ and $\sP_\lambda^{(l)}$ for the respective
eigenprojections.
If, additionally, the
inequality~{\rm(\ref{InvSquareConvergent})} holds and all the
eigenvalues $\lambda_i^{(A_1)}$ are simple, then one can
renumber the eigenprojections $\sP_\lambda^{(l)}$ in
Eq.~{\rm(\ref{sLimP})} in any way {\rm(}see
{\rm Theorem~V.4.16 of~\cite{Kato}).}
\end{remark}
\begin{theorem}\label{BasisPropertyInf}
As before, assume $\Delta_\sm=(\mu_\sm^{(1)},+\infty)$.
Also, suppose that there is a $K_B$-bounded contour
$\Gamma_l\subset D_l$ satisfying~{\rm(\ref{Best})} and such that a part
of its component $\Gamma_\sm^{l_\sm}$ coincides with the ray
$\tilde{\Delta}_\sm^{0}=[\mu_0,\ri b_0+\infty)$
where $\mu_0\in D_\sm^{l_\sm}$,
$\mu_0=a_0+\ri b_0$ with $a_0,b_0\in\R$.
Additionally, suppose
that the remaining part
\mbox{$\tilde{\Gamma}_l=\Gamma_l\setminus\tilde{\Delta}_\sm^{0}$}
of the contour $\Gamma_l$ belongs to the half-plane
$\Real\mu<a_0$, and for $\mu\in\tilde{\Delta}_\sm^{0}$
\mbox{$
\|K_B'(\mu)\|\leq\tilde{C}(1+|\Real\mu|)^{-\theta}\,,
$}
with $\tilde{C}>0$ and $\theta>1$.
Also, let the condition~{\rm(\ref{InvSquareConvergent})} be valid.
The sequence of the subspaces
$\cN_i^{(l)}=\sQ_i^{(l)}\cH_1$, \, $i=0,i_0,i_0+1,\ldots,$ \,
forms a basis of the space $\cH_1$, quadratically close
to an orthogonal one.
If, additionally, \mbox{$\mathop{\rm
dim}\sP_i^{(A_1)}\cH_1\leq n,$} for some $n\in\N$, the same for all
$i=0,i_0,i_0+1,\ldots$\,, then the union of orthonormal
vector bases of the subspaces
$\cN_i^{(l)}$, $i=0,i_0,i_0+1,\ldots,$
forms a Bari basis of the space $\cH_1$.
\end{theorem}

\vskip 1truecm

\baselineskip=12pt


\bigskip

\noindent{\small AMS Classification Numbers: Primary 47A56, 47Nxx;
                                           Secondary 47N50, 47A40.}

\begin{thebibliography}{99}
\small
\bibitem{JaffeLow}  {\sc  Jaffe, R. L.,} and  {\sc  Low, F. E.}:
        Connection between Quark--Model Eigenstates and
        Low--Energy Scattering,  Phys. Rev.
         {\bf D19} (1979), 2105--2118.
\bibitem{Simonov}  {\sc Simonov, Yu. A.}:
           Hadron--Hadron Interaction in the Compound--Bag Model,
           Yadernaya Fiz. {\bf 36} (1982), 722--731 [Russian].%
\bibitem{Okubo} {\sc Okubo, S.}: Diagonalization of Hamiltonian
         and Tamm--Dancoff equation, Progr. Theor. Phys. {\bf 12} (1954),
         603--622.%
\bibitem{GlockleMuller} {\sc Gl\"ockle, W.}, and {\sc M\"uller, L.}:
      Relativistic Theory of Interacting Particles, Phys. Rev. C
      {\bf 23} (1981), 1183--1195.%
\bibitem{KorchinShebeko} {\sc Korchin, A. Yu.,} and {\sc Shebeko, A. V.}:
     The Method of Okubo's Effective Operators and Relativistic
     Model of Nuclear Structure, Phys. At. Nucl. {\bf 56} (1993),
     1663--1671.%
\bibitem{MalyshevMinlos}
       {\sc Malyshev, V. A.,} and {\sc  Minlos, R. A.}:
       Invariant Subspaces of Clustering Operators. I.,
       J. Stat. Phys. {\bf 21} (1979), 231--242;
       Invariant Subspaces of Clustering Operators. II.,
       Comm. Math. Phys. {\bf 82} (1981), 211--226.%
\bibitem{PavlovShushkov} {\sc Pavlov, B. S.}, and
     {\sc Shushkov, A. A.}: The Theory of Extensions,
     and Null-Range Potentials with Internal Structure,
     Math. USSR Sb. {\bf 65} (1990), 147--184.%
\bibitem{McKellarMcCay}
       {\sc McKellar, B. H. J.,} and {\sc McKay C. M.}:
        Formal Scattering Theory for Energy--Dependent Potentials,
       { Aust. J.\,Phys.} {\bf 36} (1983),
        607--616.%
\bibitem{SchmidEW}
       {\sc Schmid, E. W.}:
       The Problem of Using Energy--Dependent Nucleon--Nucleon
       Potentials in Nuclear Physics,
       Helv. Phys. Acta {\bf 60} (1987), 394--397.%
\bibitem{MotSPbWorkshop} {\sc Motovilov, A. K.:}
     Potentials Appearing after the Removal of an
     Energy--Dependence and Scattering by Them, In:
     Proc. of the Intern. Workshop
    ``Mathematical aspects of the scattering
    theory and applications'',  St.~Petersburg
    University, St.~Petersburg, 1991. P.~101--108.%
\bibitem{MotRem}
    {\sc Motovilov, A. K.:} Removal of the Resolvent-like
    Energy Dependence from Interactions and Invariant
    Subspaces of a Total Hamiltonian, J.~Math. Phys.
    {\bf 36}  (1995), 6647--6664
    (LANL E-print {\tt funct-an/9606002});
    Elimination of Energy from
    Interactions Depending on It as a Resolvent,
    Theor. Math. Phys. {\bf 104} (1995), 989--1007
    (LANL E-print {\tt nucl-th/9505030}).%
\bibitem{BraunMA}
        {\sc Braun, M. A.}:
        On Relation between Quasipotential Equation
        and Schr\"odinger Equation, Teor. Mat. Fiz. {\bf 72} (1987),
        394--402 [Russian].%
\bibitem{AdL}
    {\sc Adamjan, V. M.,} and {\sc Langer, H.}: Spectral Properties
     of a Class of Operator-Valued Functions,
     J.~Operator Theory {\bf 33} (1995), 259--277.%
\bibitem{AdLMSr}
    {\sc Adamyan, V., Langer, H., Mennicken, R.,} and {\sc Saurer, J.:}
    Spectral Components of Selfadjoint Block Operator Matrices
    with Unbounded Entries,
    Math. Nachr. {\bf 178} (1996), \mbox{43--80.}%
\bibitem{MenShk}
      {\sc Mennicken, R.,} and {\sc Shkalikov, A. A.:}
      Spectral Decomposition
      of Symmetric Operator Matrices,
      Math. Nachr. {\bf 179} (1996), 259--273.%
\bibitem{ReedSimonIII}
     {\sc Reed, M.,} and {\sc Simon, B.}:
      Methods of Modern Mathematical Physics,
      III: Scattering theory,  Academic Press, N.Y., 1979.%
\bibitem{MotMathNach}
   {\sc Motovilov, A. K.:} Representations for the Three--Body T--Matrix,
   Scattering Matrices and Resolvent on Unphysical Energy Sheets,
    Math. Nachr. {\bf 187} (1997), 147--210
   (LANL E-print {\tt funct-an/9509003}).%
\bibitem{GK}
     {\sc Gohberg, I. C.,} and {\sc Krein, M. G.}:
     Introduction to the Theory of Linear Non-selfadjoint
     Operators, American Mathematical Society, Providence, 1988.%
\bibitem{Kato}
    {\sc Kato, T.} Perturbation Theory for Linear
    Operators, New York: Springer-Verlag,  1966.%
\bibitem{VirozubMatsaev} {\sc Virozub, A. I.,} and
    {\sc Matsaev, V. I.}: The Spectral Properties
    of a Certain Class of Selfadjoint Operator Functions,
    Funct. Anal. Appl. {\bf 8} (1974), 1--9.%
\bibitem{MarkusMatsaev}
    {\sc Markus, A. S.,} and {\sc Matsaev, V. I.:}
     On the Basis Property for a Certain Part of the Eigenvectors
     and Associated Vectors of a Selfadjoint Operator Pencil,
     Math. USSR Sb. {\bf 61} (1988), 289--307.%
\end{thebibliography}
\end{document}